\documentclass[english,abstracton, fleqn]{scrartcl}
\usepackage[T1]{fontenc}
\usepackage{array}
\usepackage{float}
\usepackage{bm}
\usepackage{multirow}
\usepackage{amsmath}
\usepackage{amssymb}
\usepackage{graphicx}
\usepackage{setspace}

\usepackage{authblk}
\usepackage[utf8]{inputenc}

\newcommand{\IR}{\mathbb{R}}
\newcommand{\bfn}{{\boldsymbol n}}
\newcommand{\bfv}{{\boldsymbol v}}
\newcommand{\bfw}{{\boldsymbol w}}
\newcommand{\bfu}{{\boldsymbol u}}
\newcommand{\bfx}{{\boldsymbol x}}

\newcommand{\bfI}{{\boldsymbol I}}
\newcommand{\bff}{{\boldsymbol f}}
\newcommand{\bfeps}{{\boldsymbol\varepsilon}}
\newcommand{\bfsig}{{\boldsymbol\sigma}}

\newcommand{\Ps}{{\boldsymbol P}_\Gamma}
\newcommand{\Psh}{{\boldsymbol P}_{\Gamma_h}}
\newcommand{\nablas}{\nabla_\Gamma}

\newcommand{\mcT}{\mathcal{T}}
\newcommand{\mcF}{\mathcal{F}}

\makeatletter

\providecommand{\tabularnewline}{\\}

\addtokomafont{disposition}{\rmfamily}
\usepackage{bm}

\usepackage{graphicx}
\usepackage{xcolor}

\@ifundefined{showcaptionsetup}{}{%
 \PassOptionsToPackage{caption=false}{subfig}}
\usepackage{subfig}
\makeatother

\usepackage{babel}
\begin{document}


\global\long\def\vec#1{\bm{#1}}

\global\long\def\trans{\mathsf{T}}

\title{Cut finite element modeling of linear membranes}

\author[1]{Mirza Cenanovicy\thanks{mirza.cenanovic@ju.se}}
\author[1]{Peter Hansbo\thanks{peter.hansbo@ju.se}}
\author[2]{Mats G. Larson\thanks{mats.larson@umu.se}}
\affil[1]{Department of Mechanical Engineering, Jönköping University, School of Engineering, SE-55111 Jönköping, Sweden}
\affil[2]{Department of Mathematics and Mathematical Statistics, Umeå University, SE-90187 Umeå, Sweden}

\maketitle
\begin{abstract}
We construct a cut finite element method for the membrane elasticity
problem on an embedded mesh using tangential differential
calculus. Both free membranes and membranes coupled to 3D elasticity are considered.
The discretization comes from a Galerkin method using the restriction of 3D basis funtions (linear or trilinear) to the
surface representing the membrane. In the case of coupling to 3D elasticity, we view the membrane as
giving additional stiffness contributions to the standard stiffness matrix resulting from the discretization of the three--dimensional continuum.
\end{abstract}

\section{Introduction}

In this paper we construct finite element methods for linearly elastic membranes embedded in three dimensional space meshed by tetrahedral or hexahedral elements.
These meshes do not in general align with the surface of the membrane which instead cuts through the elements. For
the modeling of the membrane problems we use the tangential differential calculus employed by Hansbo and Larson \cite{HaLa14} for meshed membranes. We extend this approach following Olshanskii, Reusken, and Grande \cite{OlReGr09} and construct a Galerkin method by using restrictions of the 3D
basis functions defined on the three--dimensional mesh to the surface. This approach can lead to severe ill conditioning, so we adapt a stabilization technique proposed by Burman, Hansbo, and Larson \cite{BuHaLa15} for the Laplace--Beltrami operator to the membrane problem.

The main application that we have in mind is the coupling of membranes to 3D elasticity. This allows for the modelling of reinforcements, such as shear strengthening and
adhesive layers. In the case of adhesives, the method can be further refined using the imperfect bonding approach of Hansbo and Hansbo \cite{HaHa04}, so that cut 3D elements
are used, allowing for relative motion of the continuum on either side of the adhesive. This extension is not explored in this paper but has been considered, for adhesives, in a discontinuous Galerkin setting in \cite{HaSa15}. Here we restrict the method to adding membrane stiffness to a continuous 3D approximation. The idea of adding stiffness from lower--dimensional structures is a classical approach, cf. Zienkiewicz \cite[Chapter 7.9]{Zienkiewicz77}, using element sides or edges as lower dimensional entities.
Letting the membranes cut through the elements in an arbitrary fashion considerably increases the practical modeling possibilities.

The paper is organized as follows: in Section 2 we introduce the membrane model problem and the finite element method for membranes and embedded membranes; in Section 3 we describe the implementation details of the method; and in Section 4 we present numerical results.
\section{The membrane model and finite element method}

\subsection{Tangential calculus}

In what follows, $\Gamma$ denotes a an oriented surface,
which is embedded in ${{\IR}}^{3}$ and equipped with exterior
normal $\bfn_\Gamma$. The boundary of $\Gamma$ consists of two parts, $\partial\Gamma_\text{N}$, where zero traction boundary conditions are assumed, and $\partial\Gamma_\text{D}$ where zeros Dirichelt boundary conditions are assumed.

We let $\rho$ denote the signed distance function
fulfilling $\nabla\rho\vert_\Gamma =\bfn_\Gamma$.

For a given function $u: \Gamma \to \IR$ we assume that there exists an extension $\bar{u}$, in some neighborhood of $\Gamma$, such that $\bar{u}\vert_\Gamma = u$.
 The the tangent gradient $\nabla_\Gamma$ on $\Gamma$ can be defined
by
\begin{equation}
\nablas u = \Ps \nabla \overline{u}
\label{eq:tangent-gradient}
\end{equation}
with $\nabla$ the ${{\IR}}^{3}$ gradient and $\Ps = \Ps(\bfx)$ the
orthogonal projection of $\IR^{3}$ onto the tangent plane of $\Gamma$ at $\bfx \in \Gamma$
given by
\begin{equation}
  \Ps = \bfI - \bfn_{\Gamma} \otimes \bfn_{\Gamma}
\end{equation}
where $\bfI$ is the identity matrix. 
The tangent gradient defined by \eqref{eq:tangent-gradient} is easily shown to be independent of the extension $\overline{u}$.
In the following, we shall consequently not make the distinction between functions on $\Gamma$ and their extensions when defining differential operators.

The surface gradient has three
components, which we shall denote by
\[
\nabla_\Sigma u=: \left(
\frac{\partial u}{\partial x^\Gamma} ,
\frac{\partial u}{\partial y^\Gamma} ,
\frac{\partial u}{\partial z^\Gamma}\right) .
\]
For a vector valued function $\bfv(\bfx)$, we define the
tangential Jacobian matrix as the transpose of the outer product of $\nabla_\Gamma$ and $\bfv$,
\[
\left(\nabla_\Gamma\otimes\bfv\right)^{\text{T}} :=\left[\begin{array}{>{\displaystyle}c>{\displaystyle}c>{\displaystyle}c}
\frac{\partial v_1}{\partial x^\Gamma} &\frac{\partial v_1}{\partial y^\Gamma} & \frac{\partial v_1}{\partial z^\Gamma} \\[3mm]
\frac{\partial v_2}{\partial x^\Gamma} &\frac{\partial v_2}{\partial y^\Gamma} & \frac{\partial v_2}{\partial z^\Gamma} \\[3mm]
\frac{\partial v_3}{\partial x^\Gamma} &\frac{\partial v_3}{\partial y^\Gamma} & \frac{\partial v_3}{\partial z^\Gamma}
\end{array}\right] ,
\]
the surface divergence $\nabla_{\Gamma}\cdot\bfv := \text{tr}\nabla_\Gamma\otimes\bfv$, and the in--plane strain tensor
\[
\bfeps_{\Gamma}(\bfu) := \Ps\bfeps(\bfu)\Ps ,\quad\text{where}\quad \bfeps(\bfu) := \frac12\left(\nabla\otimes \bfu + (\nabla \otimes\bfu)^{\rm T}\right) 
\]
is the 3D strain tensor.

\subsection{The membrane model}
We consider, following \cite{HaLa14}, the problem of finding $\bfu: \Gamma \rightarrow {{\IR}^3}$
such that
\begin{equation}
  \label{eq:LB}\begin{array}{r}
-\nabla_\Gamma\cdot \bfsig_{\Gamma}(\bfu)   = \bff \quad \text{on $\Gamma$,}\\[3mm]
\bfsig_\Gamma  =  2\mu \bfeps_\Gamma + {\lambda} \text{tr}\bfeps_\Gamma\, \Ps \quad \text{on $\Gamma$.}
\end{array}
\end{equation}
where $\bff:\Gamma\rightarrow {{\IR}^3}$ is a load per unit area and, with Young's modulus $E$ and Poisson's ratio $\nu$, 
\[
\mu := \frac{E}{2(1+\nu)},\quad {\lambda}:= \frac{E\nu}{1-\nu^2}
\]
are the Lam\'e parameters in plane stress.

Splitting the displacement into a normal part $u_{\text{N}} := \bfu\cdot\bfn_\Gamma$ and a tangential part
$\bfu_\text{T} := \bfu-u_{\text{N}} \bfn_\Gamma$,
the
corresponding weak statement takes the form: find
\[
\bfu \in V := \{\bfv  : v_\text{N} \in L_2(\Gamma )\;\text{and}\; \bfv_\text{T} \in [H^1(\Gamma)]^2, \; \bfv = {\mathbf 0}\; \text{on}\;\Gamma_\text{D}\},
\]
such that
\begin{equation}
a(\bm{u},\bm{v})=l(\bm{v}),\;\forall\bm{v}\in V,
\end{equation}
\label{membraneFEM}
where
\[
a(\bm{u},\bm{v})=(2\mu\bm{\varepsilon}_{\Gamma}(\bm{u}),\bm{\varepsilon}_{\Gamma}(\bm{v}))_{\Gamma}+(\lambda\nabla_{\Gamma}\cdot\bm{u},\nabla_{\Gamma}\cdot\bm{v})_{\Gamma}, \quad l(\bfv) = (\bff,\bfv)_\Gamma ,
\]
and 
\[
(\bfv,\bfw)_\Gamma = \int_\Gamma \bfv\cdot\bfw\,\text{d}\Gamma \quad\text{and}\quad (\bm{\varepsilon}_{\Gamma}(\bm{v}),\bm{\varepsilon}_{\Gamma}(\bm{w}))_\Gamma = \int_\Gamma \bm{\varepsilon}_{\Gamma}(\bm{v}):\bm{\varepsilon}_{\Gamma}(\bm{w})\,\text{d}\Gamma\] 
are the $L^2$ inner products.

\subsection{The cut finite element method}
\label{ssec:domain-and-fem-spaces}

Let $\widetilde{\mcT}_{h}$ be a quasi uniform mesh, with mesh parameter
$0<h\leq h_0$, into shape regular tetrahedra (hexahedra will be briefly discussed in Section \ref{numex}) of an open and bounded domain $\Omega$ in
$\IR^{3}$ completely containing $\Gamma$. 
On $\widetilde{\mcT}_h$, let $\phi$ be a continuous, piecewise
linear approximation of the distance function $\rho$
and define the discrete surface $\Gamma_h$ as the zero level set of
$\phi$; that is
\begin{equation}
\Gamma_h = \{ \bfx \in \Omega : \phi(\bfx) = 0 \}
\end{equation}

\begin{figure}

\begin{centering}
\includegraphics[width=0.45\textwidth]{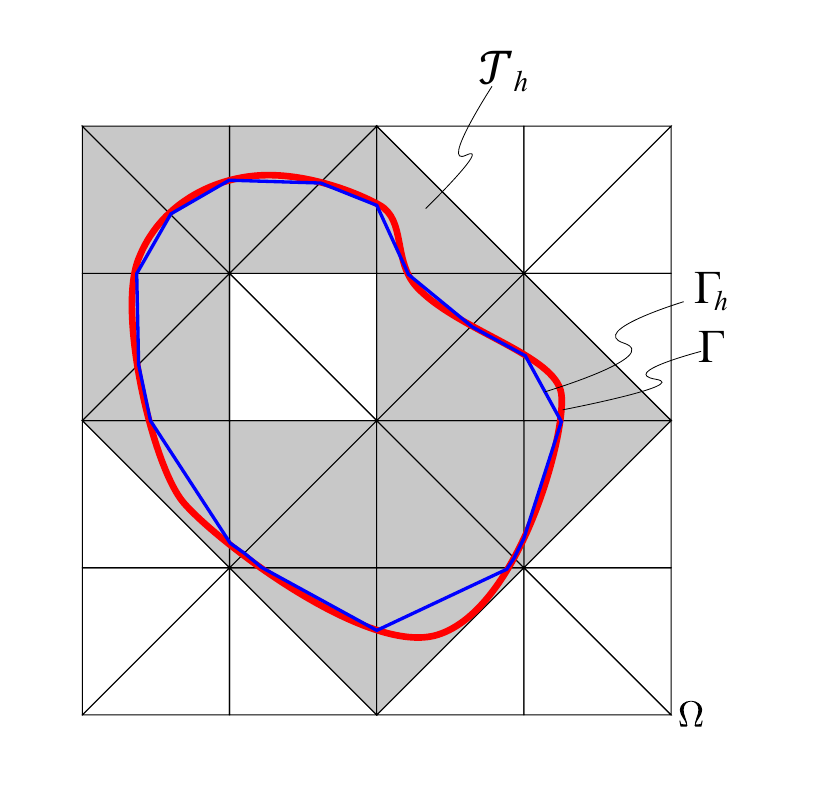}\caption{Surface domain}
\label{SurfaceDomain}
\par\end{centering}

\end{figure}

We note that $\Gamma_h$ is a polygon with flat faces and we let
$\bfn_h$ be the piecewise constant exterior unit normal to $\Gamma_h$.
For the mesh $\widetilde{\mcT}_{h}$, we define the active background
mesh by
\begin{align} 
  \mcT_h &= \{ T \in \widetilde{\mcT}_{h} : \overline{T} \cap \Gamma_h \neq \emptyset \}
  \label{eq:narrow-band-mesh}
\end{align}
cf. Figure \ref{SurfaceDomain},
and its set of interior faces by
\begin{align} 
  \mcF_h &= \{ F =  T^+ \cap T^-: T^+, T^- \in \mcT_h \}
\end{align}
The face normals $\bfn^+_F$ and $\bfn^-_F$ are then given by the unit
normal vectors which are perpendicular on $F$ and are pointing exterior
to $T^+$ and $T^-$, respectively.
%
We observe that the active background mesh $\mcT_h$ gives raise to a discrete
$h$--tubular neighborhood of $\Gamma_h$, which we denote by
\begin{align}
 \Omega_h = \cup_{T \in \mcT_h} T
\end{align}
with boundary $\partial\Omega_h$ consisting of element faces $F$ constituting the trace mesh on ${\mcT}_{h}$.
Note that for all elements
$T \in \mcT_h$ there is a neighbor $T'\in \mcT_h$ such that $T$ and
$T'$ share a face. We denote by $\partial\Omega_{h,\text{D}}$ the boundary of $\Omega_h$ intersected by $\Gamma_{h,\text{D}}$ 

Finally, let $\widetilde{V}_h$ denote the space of continuous piecewise linear (or trilinear) polynomials defined on $\widetilde{\mcT}_{h}$ and
\begin{equation}
V_h = \left\{\bfv\in [\widetilde{V}_h\vert_{\Omega_h}]^3: \; \bfv=\mathbf{0}\;\text{on}\;\partial\Omega_{h,\text{D}}\right\}
\label{eq:Vh-def}
\end{equation}
be the space of continuous piecewise linear polynomials defined
on $\mcT_h$.
The finite element method 
on $\Gamma_h$ takes the form: find $\bfu_h \in V_h$ such that
\begin{equation}\label{eq:fem}
A_h(\bfu_h,\bfv) = l_h(\bfv) \quad \forall \bfv \in V_h
\end{equation}
Here the bilinear form $A_h(\cdot, \cdot)$ is defined by
\begin{equation}
A_h(\bfv,\bfw) = a_h(\bfv,\bfw) + j_h(\bfv,\bfw) \quad \forall \bfv,\bfw \in V_h\label{eq:BilinearMembrane}
\end{equation}
with
\begin{equation}
a_h(\bfv,\bfw) = (2\mu\bm{\varepsilon}_{\Gamma_h}(\bm{v}),\bm{\varepsilon}_{\Gamma_h}(\bm{w}))_{\Gamma_h}+(\lambda\nabla_{\Gamma_h}\cdot\bm{v},\nabla_{\Gamma_h}\cdot\bm{w})_{\Gamma_h}
\end{equation}
and
\begin{equation}
j_h(v,w) = \sum_{F \in \mcF_h} \int_{F}\tau_0\left[\bfn_F^+\cdot\nabla \bfv\right] \cdot \left[\bfn_F^+\cdot\nabla \bfw\right]\, ds .
\end{equation}
Here 
$[\bfv] = \bfv^+ - \bfv^-$, where
$w(\bfx)^\pm = \lim_{t\rightarrow 0^+} w(\bfx \mp t \bfn_F^+)$, denotes
the jump of $\bfv$ across the face $F$, and $\tau_0$ is a constant of $O(1)$. 
The tangent 
gradients are defined using the normal to the discrete surface
\begin{equation}
\nabla_{\Gamma_h} v = \Psh \nabla v = ( \bfI - \bfn_h \otimes \bfn_h ) \nabla v 
\end{equation}
and the right hand side is given by
\begin{equation}
l_h(\bfv) = ( \bff, \bfv )_{\Gamma_h}
\end{equation}

\subsection{The case of embedding in a three-dimensional body}

We next consider the case of a membrane embedded in a surrounding elastic matrix. This model could be used for computation of
adhesive interfaces or reinforcements using fibers or fiber plates. The setting is very general and allows for both 2D and 1D
models to be added to the 3D model. Here we only consider adding membrane stiffness to an elastic 3D matrix,
and we then use the triangulation $\widetilde{\mcT}_{h}$ for the discretization of three--dimensional elasticity. 
In $\Omega\setminus\Gamma$ we thus assume there holds
\begin{equation}
  \label{eq:elast}
-\nabla\cdot \bfsig(\bfu)   = \bff_\Omega ,\quad
\bfsig  =  2\mu \bfeps + {\lambda_\Omega} \text{tr}\bfeps\, \bfI ,
\end{equation}
for given body force $\bff_\Omega$, where  $\lambda_\Omega := E\nu/((1+\nu)(1-2\nu))$.
We assume for simplicity of presentation that $\bfu={\mathbf 0}$ on $\partial\Omega_\text{D}$, a part of the boundary which is assumed to include $\partial\Omega_{h,\text{D}}$, 
and that the traction is zero on the rest of the boundary. Our finite element method in the bulk is then based on the finite element space
\begin{equation}
W_h = \left\{\bfv\in [\widetilde{V}_h]^3: \; \bfv=\mathbf{0}\;\text{on}\;\partial\Omega_{\text{D}}\right\},
\label{eq:Wh-def}
\end{equation}
and we seek $\bfu_h\in W_h$ such that
\begin{equation}\label{eq:femfull}
a_{\Omega}(\bfu_h,\bfv)+
a_h(\bfu_h,\bfv) =l_{\Omega}(\bfv)+ l_h(\bfv) \quad \forall \bfv \in W_h
\end{equation}
where
\[
a_{\Omega}(\bfu_h,\bfv):=(2\mu\bm{\varepsilon}(\bm{v}),\bm{\varepsilon}(\bm{w}))_{\Omega}+(\lambda\nabla\cdot\bm{v},\nabla\cdot\bm{w})_{\Omega}, \quad l_{\Omega}(\bfv):=(\bff_\Omega,\bfv)_{\Omega}.
\]
The FEM (\ref{eq:femfull}) thus takes into account both the stiffness from the bulk and from the membrane. The bulk stiffness matrix is here established independently of the position of the membrane which allows for rapid repositioning of the membrane; this is beneficial for example for the purpose of optimizing the membrane location.
We remark that when the membrane is embedded in a three--dimensional mesh used for elasticity in all of $\Omega$, then we can drop the stabilization term (or set $\tau_0=0$) since
the three--dimensional stiffness matrix gives stability to the embedded membrane.

In case of cohesive zone modelling at the membrane we also need to allow for independent relative motion of the bulk on either side of the membrane, which requires cutting the bulk elements following \cite{HaHa04}. In such a case we also need to reintroduce the stabilization term $j_h(\cdot,\cdot)$. We will return to this problem in future work.

\section{Implementation\label{numex}}

This  Section describes the implementational aspects of the embedded
membrane model and provides an algorithm of the implementation. 

\begin{figure}

\begin{centering}
\includegraphics[width=0.45\textwidth]{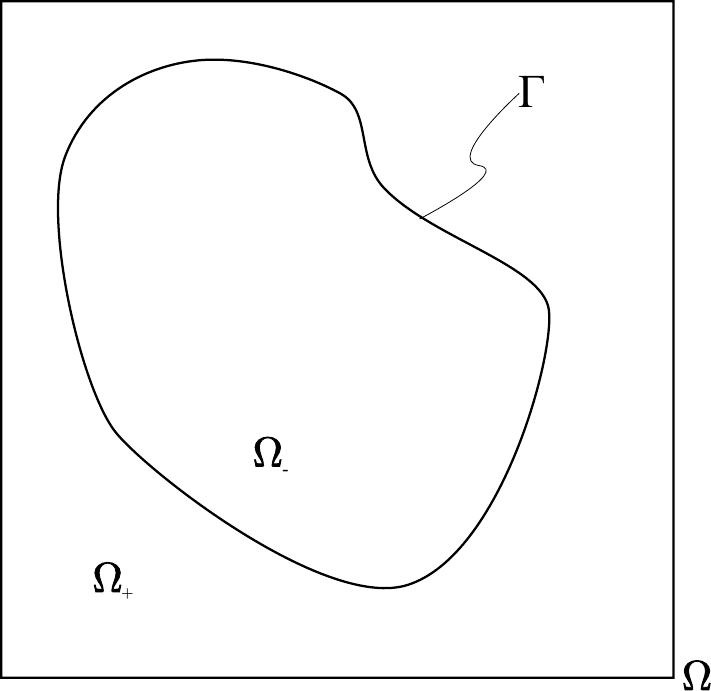}\caption{2D representation of the problem domain}
\label{ProblemDomain}
\par\end{centering}

\end{figure}

\subsection{Algorithm\label{sub:Algorithm-Membrane}}
\begin{figure}
\begin{centering}
\subfloat[2D case]{\begin{centering}
\includegraphics[width=0.3\columnwidth]{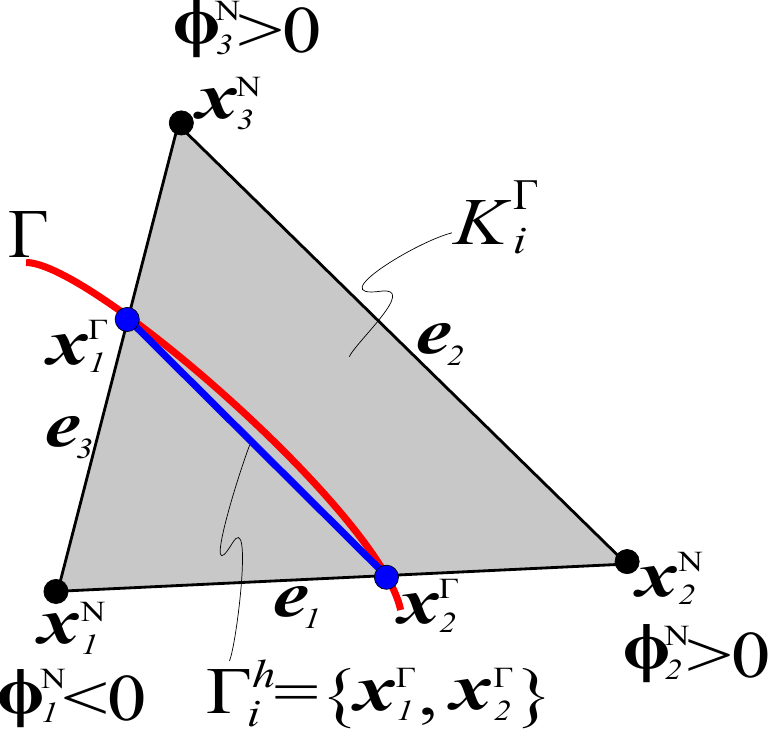}
\par\end{centering}

}\hspace{0.02\textwidth}\subfloat[Surface element normal]{\centering{}\includegraphics[width=0.3\textwidth]{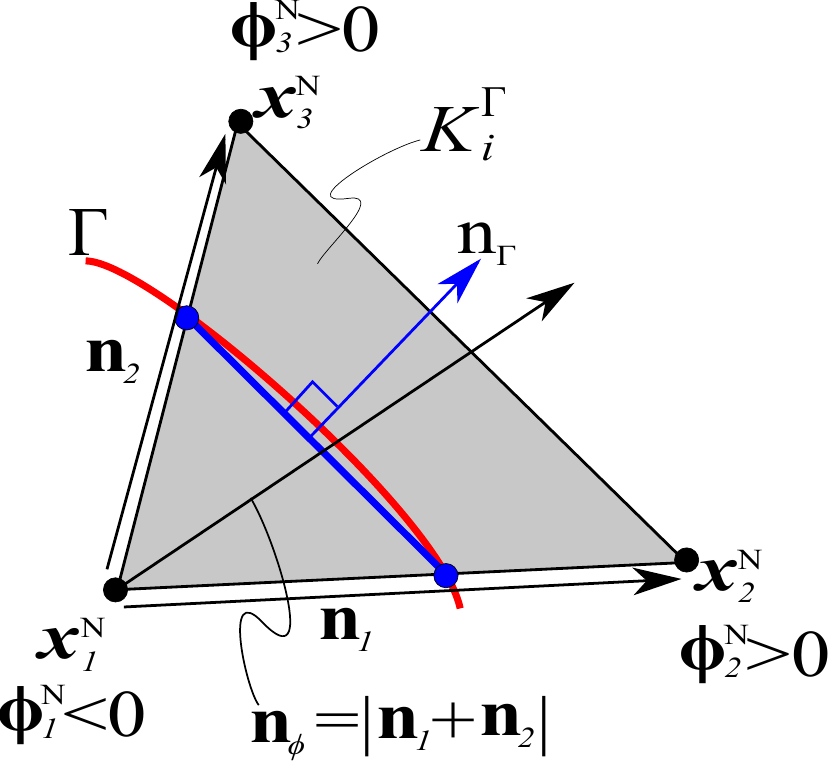}\label{SurfaceElementNormal}}\hspace{0.02\textwidth}\subfloat[3D case]{\centering{}\includegraphics[width=0.3\columnwidth]{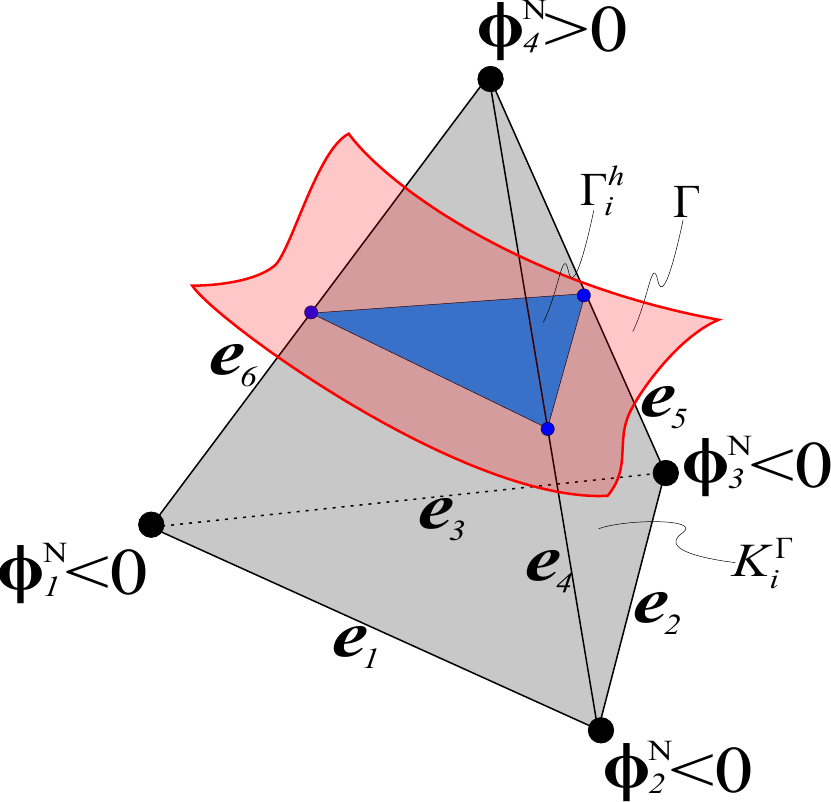}}\caption{Surface element $\Gamma_{h}^{i}$ and parent element $K_{\Gamma}^{i}$
in 2D and 3D}
\label{SurfaceElement}
\par\end{centering}

\end{figure}

\begin{enumerate}
\item Construct a mesh $\tilde{\mathcal{T}}$ in $\mathbb{R}^{d}$ on the
domain $\Omega$ in which the implicit surface $\Gamma$ will be embedded.
Let $\mathbf{x}_{N}$ denote the vector of coordinates in $\tilde{\mathcal{T}}$.
\item Construct the level set function $\rho(\bm{x})$ either analytically
or by the use of surface reconstruction, see  Section(\ref{sub:Implicit-surface-construction})
for details. \label{enu:Construct-the-level-set-function}
\item Discretize the distance function $\phi=\rho(\mathbf{x}_{N})$ by evaluating
$\rho$ in the nodes of the complete underlying mesh $\tilde{\mathcal{T}}$.
\item Find the indices to the elements in the background mesh $\mathcal{T}_{h}$,
by using the discrete distance values $(\phi>0\textnormal{ and }\phi<0)$
in some nodes of element $T_{i}$, see Figure \ref{SurfaceElement}.
\item Extract the following surface quantities. For every element $T_{i}\in\mathcal{T}_{h}$

\begin{enumerate}
\item Find the zero surface points of the polygons $\Gamma_{h}$ by looping
over all elements $T_{i}\in\mathcal{T}_{h}$ and interpolating the
discrete signed distance function $\phi$ linearly for every element
edge that has a difference in function value (see Figure \ref{SurfaceElement}
and Section \ref{sub:Zero-level-surface} for details). For simplicity
the polygon $\Gamma_{h,i}$ of element $T_{i}$ is split into triangular
elements $\hat{T}$. 
\item Compute the face normal ${\bf n}_{f}$ of each triangular element,
which is used to compute the Jacobian for the basis functions of element
$T_{i}$. Note that care must be taken when definining the face normal,
such that the orientation of the normal field becomes unidirectional.
The element normal ${\bf n}_{\phi}$ can be used to orient ${\bf n}_{f}$
in the same direction.
\end{enumerate}
\item Compute the displacement field $\bm{u}_{h}$ on $\mathcal{T}_{h}$
by solving linear system that results from the bilinear equation (\ref{eq:BilinearMembrane}).
\item Interpolate the solution field $\mathbf{u}_{h}$ to $\mathbf{u}_{\Gamma_{h}}$
using the basis functions of the elements in $\mathcal{T}_{h}$.
\end{enumerate}

\subsection{Implementation details}

\subsubsection{Implicit surface construction\label{sub:Implicit-surface-construction}}

There are a number of ways to generate an implicit surface for analysis.
An implicit surface can be approximated from a CAD surface using surface
reconstruction techniques, see Belytschko et al \cite{Belytschko2003}.
Another way is to use analytical implicit surfaces descriptions, see
Burman et al \cite{Burman2014} for an overview. In this paper we
use the following analytical surface descriptions.

\begin{figure}
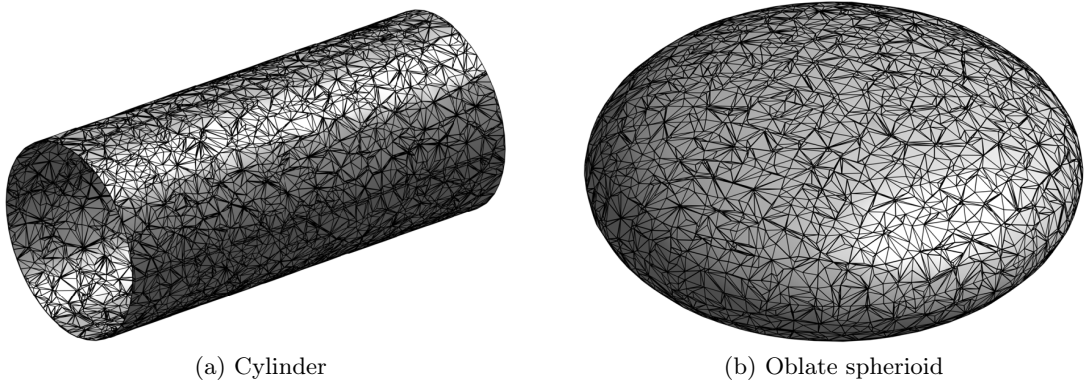

\begin{centering}
\subfloat[Cylinder]{\begin{centering}
\includegraphics[width=0.45\textwidth]{CylinderSurface}\label{CylinderSurface}
\par\end{centering}

}\hspace{0.05\textwidth}\subfloat[Oblate spherioid]{\centering{}\includegraphics[width=0.45\textwidth]{OblateSurface}\label{OblateSurface}}\caption{Implicit surfaces}
\label{ImplicitSurfaces}
\par\end{centering}

\end{figure}

Cylinder function:

\begin{equation}
\rho(x,y,z)=\sqrt{(x-x_{c})^{2}+(y-y_{c})^{2}}-r\label{eq:CylinderImplicitFunction}
\end{equation}

where $[x_{c},y_{c}]$ is the center of the cylinder. See Figure \ref{CylinderSurface}.

Oblate spherioid:

\begin{equation}
\rho(x,y,z)=x^{2}+y^{2}+(2z)^{2}-1
\end{equation}

See Figure \ref{OblateSurface}.

\subsubsection{Zero level surface approximation\label{sub:Zero-level-surface}}

\begin{figure}
\begin{centering}
\subfloat[Triangular ]{\centering{}\includegraphics[width=0.22\textwidth]{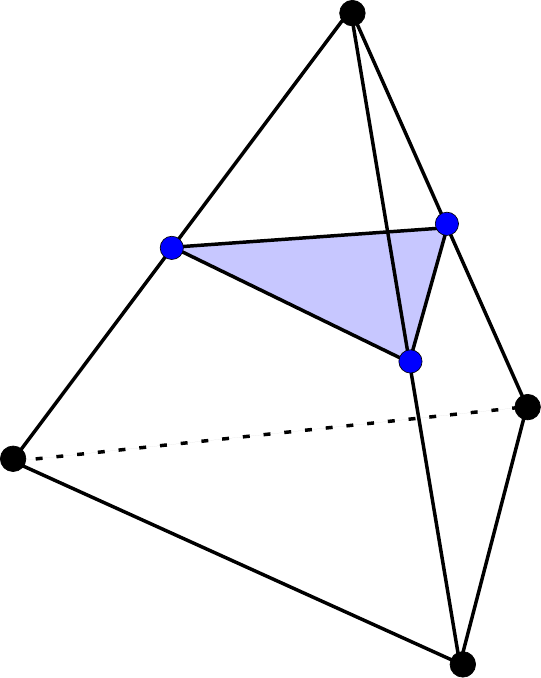}}\hspace{0.05\textwidth}\subfloat[Quadliteral ]{\centering{}\includegraphics[width=0.22\textwidth]{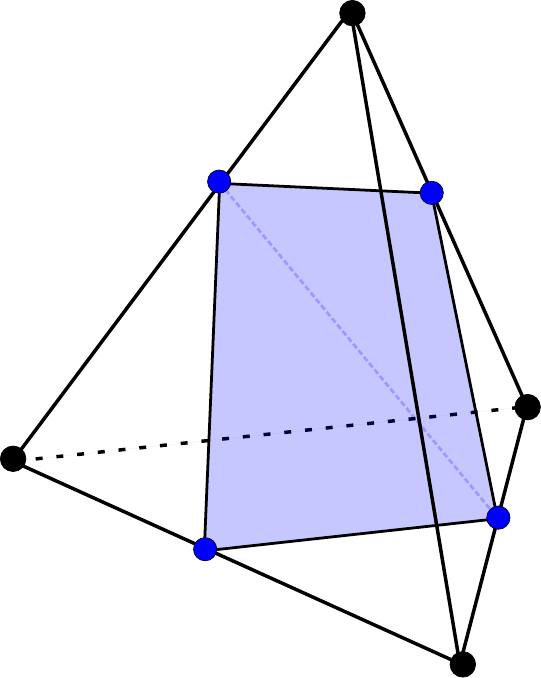}}
\par\end{centering}

\caption{Tetrahedral cut cases}
\label{TetCutCases}
\end{figure}

\begin{figure}
\begin{centering}
\subfloat[Triangular ]{\begin{centering}
\includegraphics[width=0.22\textwidth]{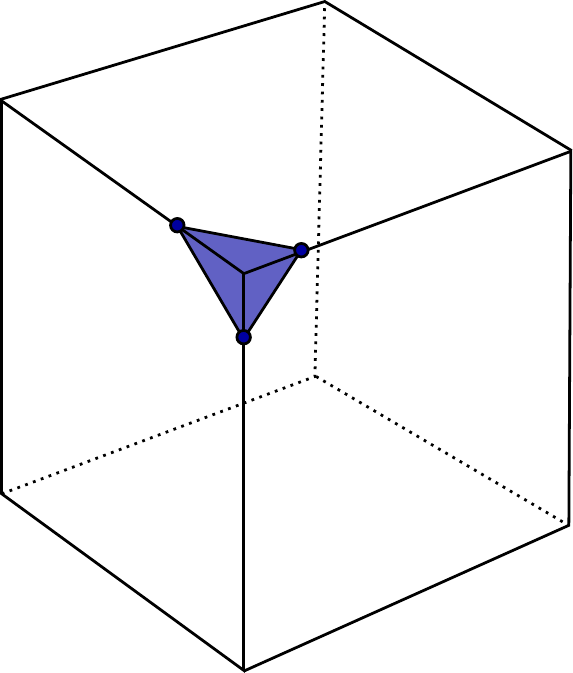}
\par\end{centering}

}\hspace{0.01\textwidth}\subfloat[Quadliteral ]{\centering{}\includegraphics[width=0.22\textwidth]{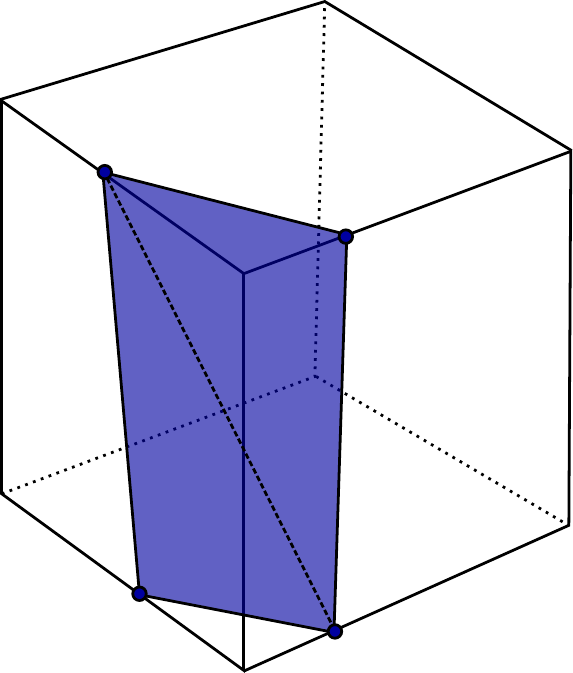}}\hspace{0.01\textwidth}\subfloat[Pentagon ]{\centering{}\includegraphics[width=0.22\textwidth]{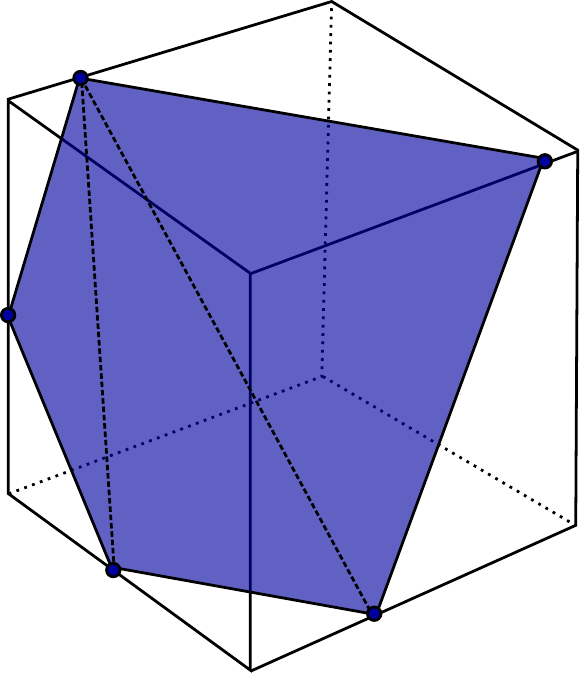}}\hspace{0.01\textwidth}\subfloat[Hexagon ]{\centering{}\includegraphics[width=0.22\textwidth]{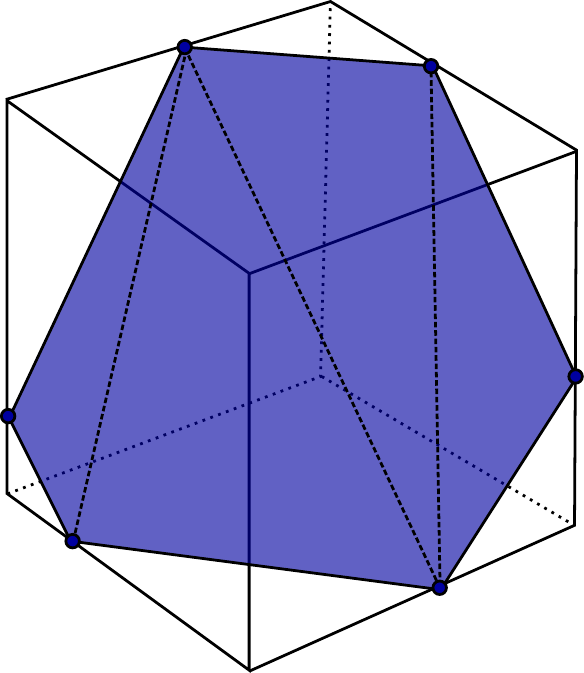}}
\par\end{centering}

\caption{Hexahedral cut cases}
\label{HexCutCases}
\end{figure}

The overal procedure is to use linear interpolation on the discrete
interface values for each element edge in order to find the zero level
surface point, see Figure \ref{SurfaceElement}. 
\begin{enumerate}
\item For every element $T_{i}\in\mathcal{T}_{h}$ loop over all edges $e_{i}$.

\begin{enumerate}
\item For every edge $e_{i}$ check the sign of the two discrete function
values $\phi|_{e_{i}}$ to determine if the edge is cut.
\item Linearly interpolate the cut point ${\bf x}_{\Gamma,i}$ along the
edge $e_{i}$ using the two vertex coordinates ${\bf x}_{e_{i}}^{m}$
and ${\bf x}_{e_{i}}^{n}$, at nodes $m$ and $n$ (endpoints of $e_{i}$)
and the function values $\phi|_{e_{i}}=\{\phi({\bf x}_{e_{i}}^{m}),\phi({\bf x}_{e_{i}}^{n})\}$.
\item Let ${\bf x}_{e_{i},1}={\bf x}_{e_{i}}|_{\phi|_{e_{i}}>0}$ (the coordinate
corresponding to the highest value of $\phi$) and ${\bf x}_{e_{i},0}={\bf x}_{e_{i}}|_{\phi|_{e_{i}}<0}$
and compute the vector ${\bf {\bf n}}_{i}={\bf x}_{e_{i},1}-{\bf x}_{e_{i},0}$.
See Figure \ref{SurfaceElementNormal}
\end{enumerate}
\item Compute the element vector $\bm{n}_{\phi}=\sum\bm{n}_{i}$.
Note that $\bm{n}_{\phi}$ points in the general direction of $\nabla\phi$ and
is only used to determine the orientation of the face normals.
\item Depending on the number of nodes in element $K_{i}^{\Gamma}$ and
the orientation of the cut, several cut cases must be considered,
see Figure \ref{TetCutCases} for tetrahedral element and Figure \ref{HexCutCases}
for hexahedral. A tessellation into triangles is done for all cases. 
\item In order to do the tessellation into triangles from an arbitrary polygon,
a rotation from $\mathbb{R}^{3}$ into $\mathbb{R}^{2}$ was performed
and then a 2D convex hull algorithm was applied.
\end{enumerate}

\subsection{Membrane embedded in elastic material}

In this  Section we demonstrate one particular application of the elastic
membrane. Consider a set of membrane surfaces embedded into an elastic
material body. We let the embedded membrane stiffen the solution by
adding stiffness from the membrane solution to the bulk solution.
This is easily done since the membrane shares the same degrees of
freedom as the bulk.

\subsubsection{Algorithm}

The algorithm described here is similar to the one in  Section \ref{sub:Algorithm-Membrane}.
We allow for several membranes (with the possibility of different
material properties) to stiffen the bulk.
\begin{enumerate}
\item Construct a mesh $\tilde{\mathcal{T}}$ in $\mathbb{R}^{d}$ on the
elastic domain $\Omega$ in which the implicit surface $\Gamma$ will
be embedded. Let $\mathbf{x}_{N}$ denote the vector of coordinates
in $\tilde{\mathcal{T}}$.
\item Create a set of surfaces functions $\{\rho(\bm{x})\}$ in the same
way as in the previous algorithm.
\item Discretize the distance functions $\{\phi\}=\{\rho(\mathbf{x}_{N})\}$
by evaluating all functions in the set $\{\rho\}$ in the nodes of
the complete underlying mesh $\mathcal{\tilde{\mathcal{T}}}$. 
\item Find the sets $\{\mathcal{T}_{h}\}$ that are intersected by the surfaces
such that for each $\{\phi\}_{i}$ there exists a corresponding set
of cut elements $\{\mathcal{T}_{h}\}_{i}$. 
\item Follow the same approach as described in the previous algorithm to
extract the zero surface information for each $\{\mathcal{T}_{h}\}_{i}$.
\item Create a discrete system of equations for the bulk elasticity problem.
While assembling the bulk stiffness matrix, for each element that
is cut by the a membrane surface, compute the membrane element stiffness
and add it to the global bulk stiffness matrix. Note that no stabilization
is needed in this case since the surrounding elements create a well
conditioned stiffnessmatrix.
\item Solve the linear system.\end{enumerate}

\section{Numerical examples}

In this Section we show convergence comparison on some geometries
for which we can compute the solutions analytically. We compare the
convergence rates of this approach with a triangulated surface. Numerical
examples are given for both tetrahedral and hexahedral elements.

The meshsize is defined by:

\[
h:=\dfrac{1}{\sqrt[3]{NNO}}
\]
where $NNO$ denotes the total number of nodes on the underlying 3D
mesh $\widetilde{\mcT}_{h}$, which is uniformly refined.

\subsection{Pulling a cylinder}

Comparing this approach to the approach previously done by Hansbo
and Larson \cite{HaLa14}, we consider a cylindrical shell of radius
$r$ and thickness $t$, with open ends at $x=0$ and $x=L$ and with
fixed axial displacements at $x=0$ and radial at $x=L$, carrying
a axial surface load per unit area
\[
f(x,y,z)=\frac{F}{2\pi r}\frac{x}{L^{2}},
\]
where $F$ has the unit of force. The resulting axial stress is
\[
\sigma=\frac{F(1-(x/L)^{2})}{4\pi rt}.
\]

We consider the same example as was used in \cite{HaLa14} and choose
$r=1$, $L=4$, $E=100$, $\nu=1/2$, $t=10^{-2}$ and $F=1$.

\begin{figure}[H]
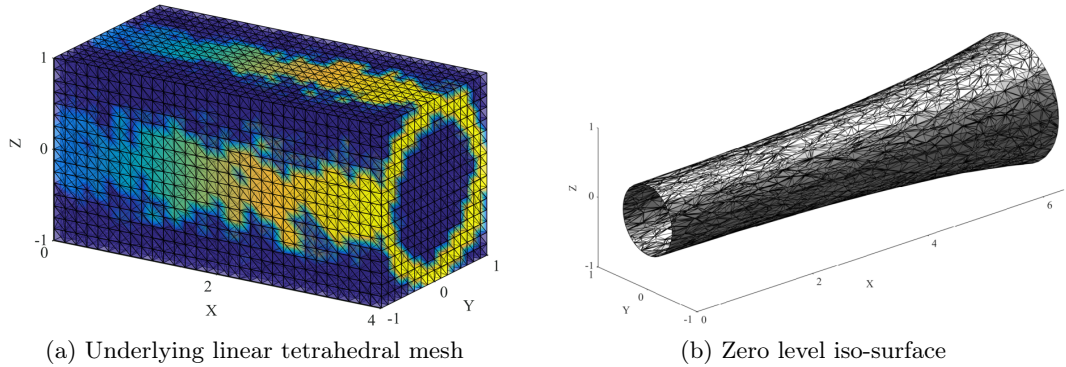

\centering{}\subfloat[Underlying linear tetrahedral mesh]{\centering{}\includegraphics[width=0.45\textwidth]{CylinderUnderlyingSolutionTet}}\hspace{0.05\textwidth}\subfloat[Zero level iso-surface]{\centering{}\includegraphics[width=0.45\textwidth]{CylinderSolutionTet}}\caption{Displacement field on a tetrahedral mesh (x10 exaggerated)}
\label{PullCylinderTetSolution}
\end{figure}

\begin{figure}[H]
\centering{}\includegraphics[width=0.75\textwidth]{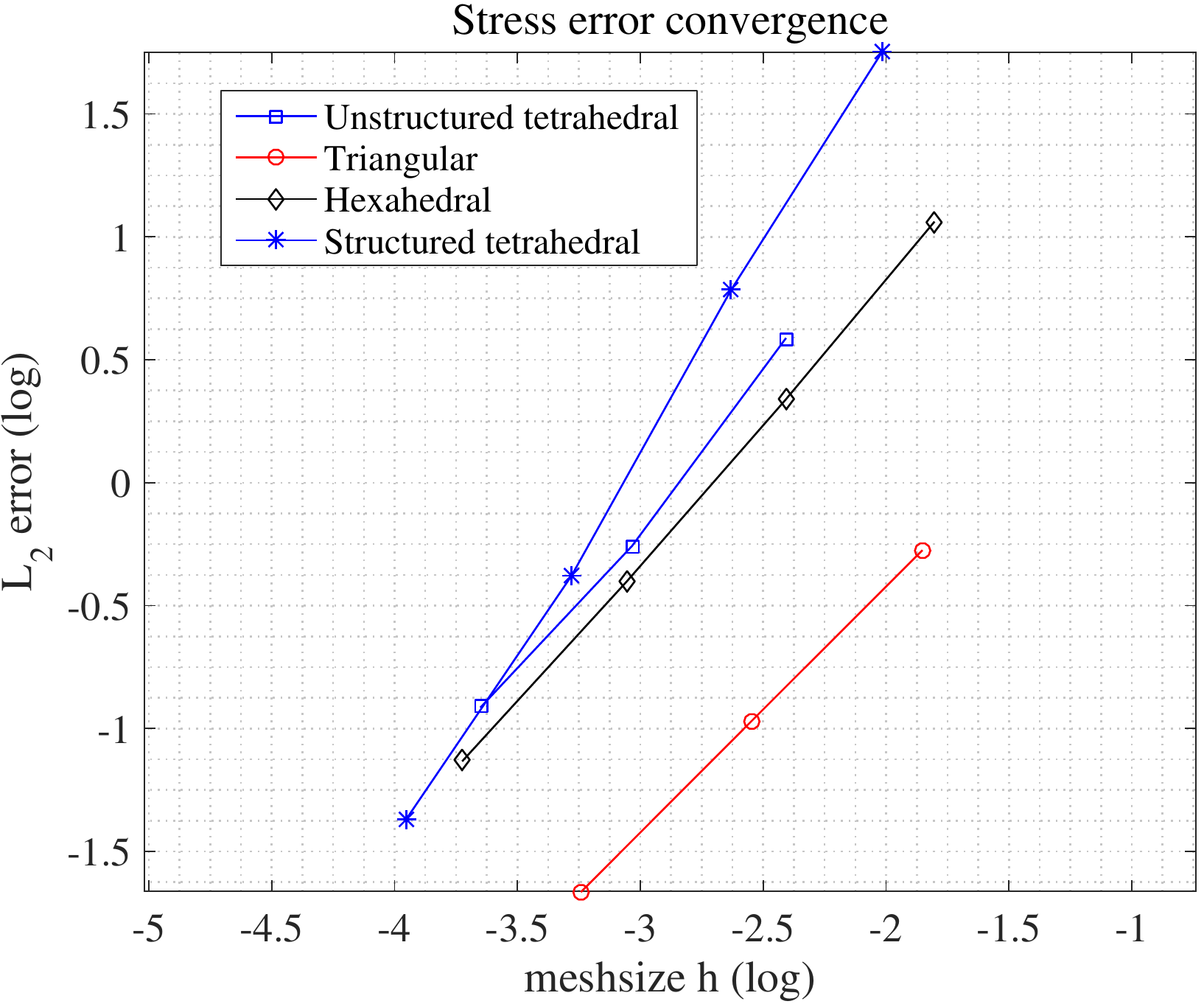}\caption{Stress convergence comparison for the cylinder.}
\label{PullCylinderConvergence}
\end{figure}

In Figure \ref{PullCylinderTetSolution} we show the solution (exaggerated
10 times) on a given tetrahedral mesh and in Figure \ref{PullCylinderHexSolution} the corresponding
solution on a hexahedral mesh with the same mesh size. In Figure \ref{PullCylinderConvergence}
we show the $L_{2}(\Sigma)$ error in stress, $||\bm{\sigma}-\bm{\sigma}_{h}||$,
where $\bm{\sigma}:=\bm{\sigma}_{\Gamma}(\bm{u})$ and $\bm{\sigma}_{h}:=\bm{\sigma}_{\Gamma}(\bm{u}_{h})$.
See table \ref{CylinderConvergenceTable} for convergence rates.

\begin{figure}[H]
\begin{centering}
\subfloat[Displement field on the underlying mesh. (Exaggerated lengths)]{\centering{}\includegraphics[width=0.45\columnwidth]{MeshDependency1}}\subfloat[Displacement field on the extracted linear zero level iso-surface.]{\begin{centering}
\includegraphics[width=0.45\columnwidth]{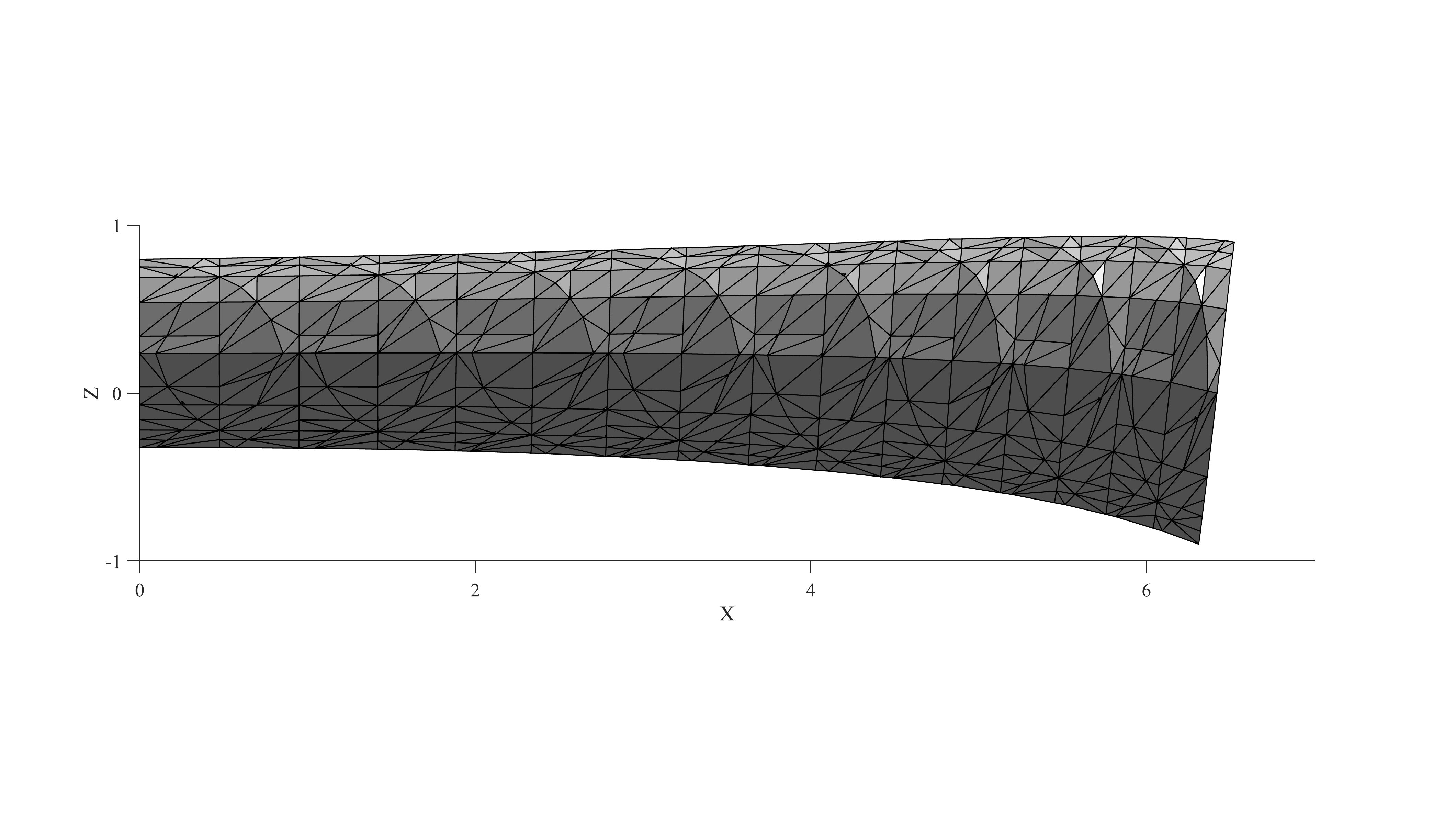}
\par\end{centering}

}
\par\end{centering}

\centering{}\caption{Mesh dependent errors on low resolution linear tetrahedral mesh.}
\label{MeshDependErrorTet}
\end{figure}

\begin{figure}[H]
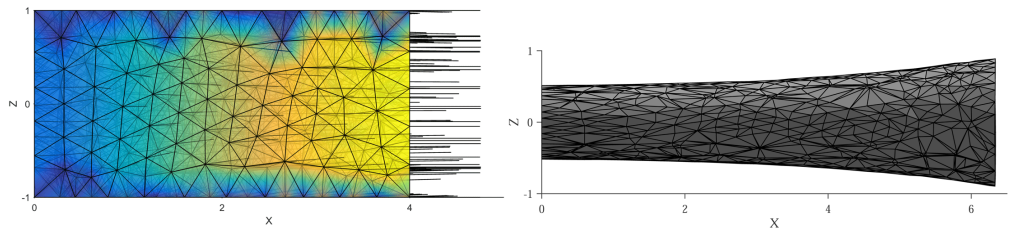

\centering{}\subfloat[Underlying displacement field]{\centering{}\includegraphics[width=0.45\textwidth]{UnstructeredUnderlying}}\subfloat[Surface displacement field]{\centering{}\includegraphics[width=0.45\textwidth]{Unstructered}}\caption{Displacement fields on unstructered tetrahedral mesh}
\label{MeshDepErrorUnstructTet}
\end{figure}

\begin{figure}[H]
\centering{}\includegraphics[width=0.45\textwidth]{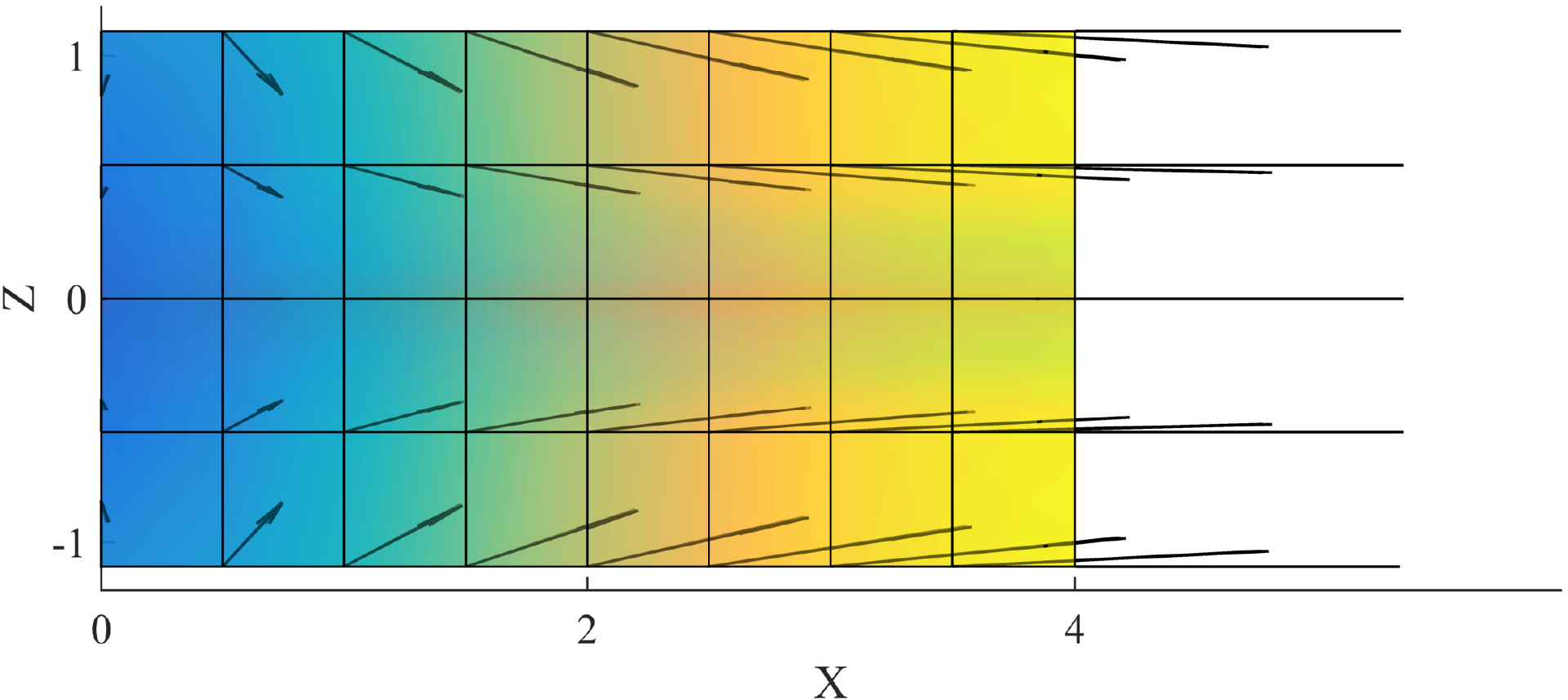}\caption{Displacement field on the underlying tri-linear hexahedral mesh mesh.
(x10 exaggerated)}
\label{UnderlyingSolutionHex}
\end{figure}

Mesh dependent errors occur in a structured tetrahedral mesh case,
see  Figure \ref{MeshDependErrorTet}, the error becomes less prominent
with a finer mesh. This error can be avoided by using an unstructured
tetrahedral mesh, see Figure \ref{MeshDepErrorUnstructTet}. This
error is not found in the hexahedral case, see Figure \ref{UnderlyingSolutionHex}.

\begin{figure}[H]
\centering{}\subfloat[Underlying linear hexahedral mesh]{\centering{}\includegraphics[width=0.45\textwidth]{CylinderUnderlyingSolutionHex}}\hspace{0.05\textwidth}\subfloat[Zero level iso-surface]{\centering{}\includegraphics[width=0.45\textwidth]{CylinderSolutionHex}}\caption{Displacement field on a hexahedral mesh (x10 exaggerated)}
\label{PullCylinderHexSolution}
\end{figure}

\begin{table}[H]
\begin{onehalfspace}
\noindent \centering{}%
\begin{tabular}{>{\raggedright}p{0.15\textwidth}>{\centering}p{0.15\textwidth}>{\centering}p{0.25\textwidth}>{\centering}p{0.25\textwidth}}
\hline 
Element type & Mesh size & $||\sigma-\sigma_{h}||$ & Rate\tabularnewline
\hline 
\hline 
\multirow{3}{0.15\textwidth}{Tetrahedral unstructured} & 0.0899 & 1.7959 & -\tabularnewline
\cline{2-4} 
 & 0.0481 & 0.7733 & 1.3498\tabularnewline
\cline{2-4} 
 & 0.0260 & 0.4017 & 1.0650\tabularnewline
\hline 
\multirow{4}{0.15\textwidth}{Tetrahedral structured} & 0.1330 & 5.7545 & -\tabularnewline
\cline{2-4} 
 & 0.0721 & 2.1916 & 1.5749\tabularnewline
\cline{2-4} 
 & 0.0376 & 0.6846 & 1.7903\tabularnewline
\cline{2-4} 
 & 0.0192 & 0.2550 & 1.4722\tabularnewline
\hline 
\multirow{4}{0.15\textwidth}{Hexahedral} & 0.1644 & 2.8832 & -\tabularnewline
\cline{2-4} 
 & 0.0899 & 1.4008 & 1.1954\tabularnewline
\cline{2-4} 
 & 0.0472 & 0.6698 & 1.1438\tabularnewline
\cline{2-4} 
 & 0.0242 & 0.3228 & 1.0924\tabularnewline
\hline 
\multirow{3}{0.15\textwidth}{Original} & 0.1565 & 0.7580 & -\tabularnewline
\cline{2-4} 
 & 0.0783 & 0.3793 & 0.9998\tabularnewline
\cline{2-4} 
 & 0.0392 & 0.1897 & 1.0015\tabularnewline
\cline{2-4} 
\end{tabular}\caption{Convergence of the cylinder}
\label{CylinderConvergenceTable}\end{onehalfspace}
\end{table}

\subsection{Pulling an oblate spherioid}

\begin{figure}[H]
\centering{}\includegraphics[width=0.5\textwidth]{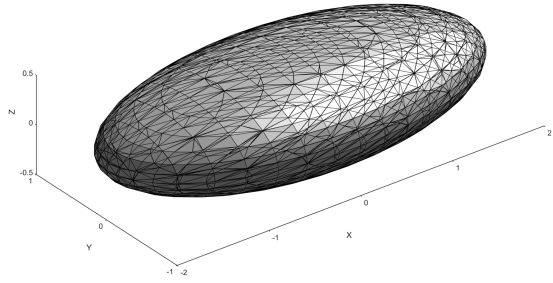}\caption{Displacement field on the extracted linear zero level iso-surface
from a linear tetrahedral mesh.(x10 exaggerated)}
\label{DisplFieldOblateTet}
\end{figure}

\begin{figure}[H]
\centering{}\includegraphics[width=0.75\textwidth]{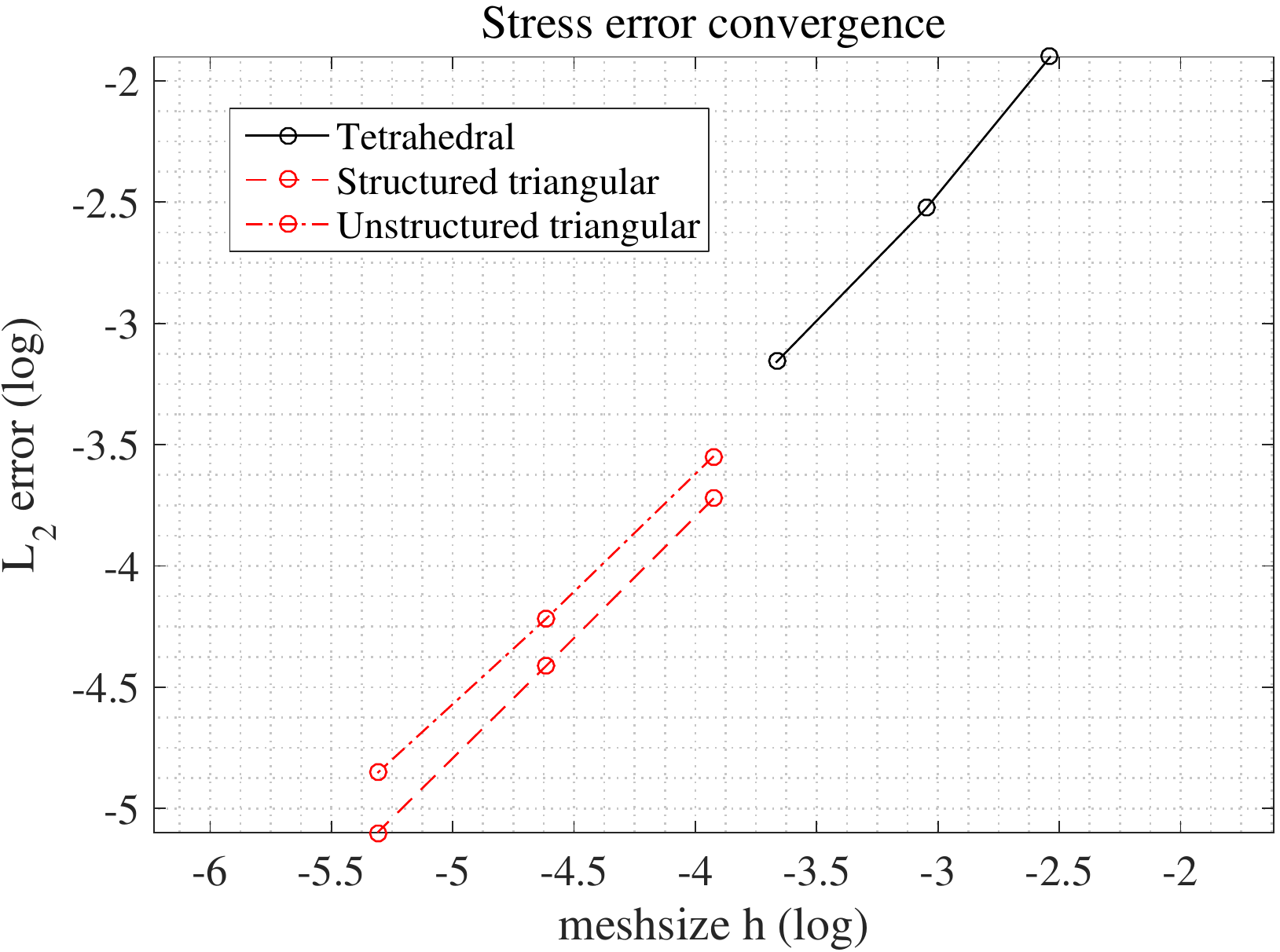}\caption{Convergence of stress}
\label{OblateStressConvergence}
\end{figure}

Again we use the same example as in \cite{HaLa14} an set the exact
solution to be $\bm{u}=(x,0,0)$ and compute the stress and then the corrsponding load from (\ref{eq:LB}). We set the parameters $\textrm{E}=1$, $\nu=1/2$, and $t=1$. The computed displacement
field can be seen in Figure \ref{DisplFieldOblateTet}. Compared to
the previous work in \cite{HaLa14}, the superparametric stabilization
method is not needed in this approach since we already stabilize using
the term $j_h(\cdot,\cdot)$. A stress convergence comparison can be seen in Figure
\ref{OblateStressConvergence}. The convergence rates can be seen
in Table \ref{ConvergenceOblateTable}. As one can see from the number
of degree of freedoms and the corresponding meshsize in each approach
(table \ref{ConvergenceOblateTable}) , the approximation of the embedded
approach is dependent on the surface curvature or the surface to volume
ratio. If the ratio is small then the surface is cutting fewer elements
and the resolution of the background mesh can be kept coarse but still
generate a good approximation.

\begin{table}[H]
\begin{onehalfspace}
\noindent \centering{}%
\begin{tabular}{>{\raggedright}p{0.15\textwidth}>{\centering}p{0.15\textwidth}>{\centering}p{0.15\textwidth}>{\centering}p{0.25\textwidth}>{\centering}p{0.25\textwidth}}
\hline 
Element type & Mesh size & Number of DOFs & $||\sigma-\sigma_{h}||$ & Rate\tabularnewline
\hline 
\hline 
\multirow{3}{0.15\textwidth}{Tetrahedral} & 0.0790 & 2414 & 0.1493 & -\tabularnewline
\cline{2-5} 
 & 0.0474 & 7218 & 0.0800 & 1.2222\tabularnewline
\cline{2-5} 
 & 0.256 & 26961 & 0.0425 & 1.0261\tabularnewline
\hline 
\multirow{3}{0.15\textwidth}{Structured, triangular} & 0.0198 & 2562 & 0.0242 & -\tabularnewline
\cline{2-5} 
 & 0.0099 & 10242 & 0.0121 & 1.0004\tabularnewline
\cline{2-5} 
 & 0.0049 & 40962 & 0.0061 & 0.9882\tabularnewline
\hline 
\multirow{3}{0.15\textwidth}{Unstructured, triangular} & 0.0198 & 2562 & 0.0288 & -\tabularnewline
\cline{2-5} 
 & 0.0099 & 10242 & 0.0147 & 0.9707\tabularnewline
\cline{2-5} 
 & 0.0049 & 40962 & 0.0078 & 0.9144\tabularnewline
\cline{2-5} 
\end{tabular}\caption{Convergence of the oblate}
\label{ConvergenceOblateTable}\end{onehalfspace}
\end{table}

\subsection{Membrane embedded in elastic material }

\begin{figure}[H]
\begin{centering}
\includegraphics[width=0.6\textwidth]{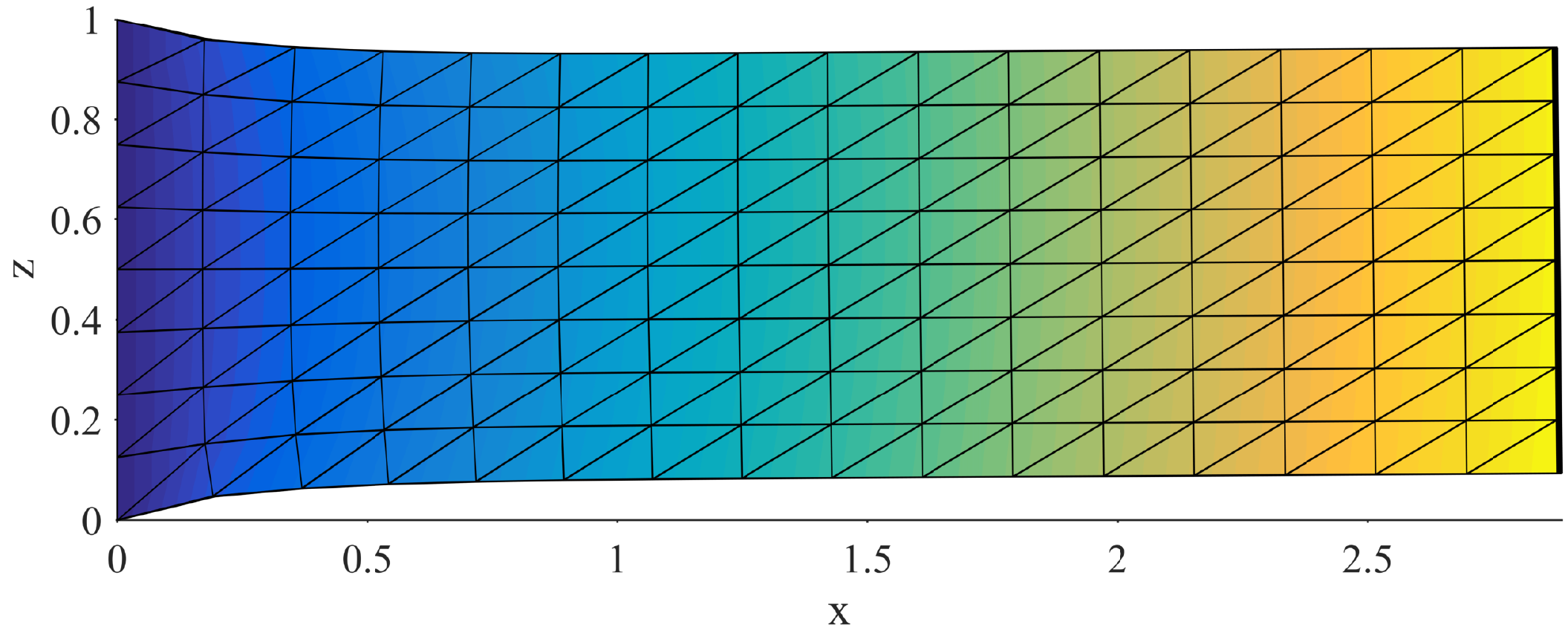}
\par\end{centering}

\caption{Elastic beam without embedded elastic membrane}
\label{ElasticBulkNoMembrane}
\end{figure}

A rectangular box $(0\leq x\leq2$ and $0\leq y,z\leq1$) has a surface
load, $f$=1, applied in the positive x direction at $x=2$. The bulk
material has the following properties, $\mathrm{E}=100,$ $\nu=0.5$.
The membrane has $\mathrm{E}=1000,$ $\nu=0.5$ and $t=0.01$. Figure
\ref{ElasticBulkNoMembrane} shows a displacement plot (40x exaggerated)
of the elastic bulk material without any embedded membrane. Figure
\ref{EmbeddedMaterial1} shows the same displacement plot but with
8 embedded membrane surfaces. The surfaces are visualized in Figure
\ref{EmbeddedMaterial2}. 

\begin{figure}[H]
\begin{centering}
\includegraphics[width=0.6\textwidth]{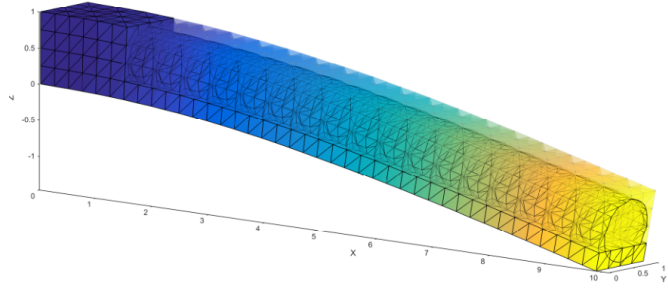}
\par\end{centering}

\caption{Elastic beam in bending with embedded elastic membrane}
\label{ElasticBulkBendingWithStiffening}
\end{figure}

Finally, in Figures \ref{ElasticBulkBendingNoStiffening} and \ref{ElasticBulkBendingWithStiffening}
we how the effect on displacements resulting from inserting a circular membrane into a beam in bending. The material properties are the same as in the previous example.
We note the marked increase in bending stiffness resulting from the added membrane stiffness.

\begin{figure}[H]
\begin{centering}
\includegraphics[width=0.6\textwidth]{stiffening-8Membrane}
\par\end{centering}

\caption{Elastic beam with embedded elastic membrane}
\label{EmbeddedMaterial1}
\end{figure}

\begin{figure}[H]
\begin{centering}
\includegraphics[width=0.6\textwidth]{stiffening-8MembraneCut}
\par\end{centering}

\caption{Elastic beam with embedded elastic membrane}
\label{EmbeddedMaterial2}
\end{figure}

\begin{figure}[H]
\begin{centering}
\includegraphics[width=0.6\textwidth]{BendingElasticBulk}
\par\end{centering}

\caption{Elastic beam in bending without embedded elastic membrane}
\label{ElasticBulkBendingNoStiffening}
\end{figure}

\section{Concluding remarks}

In this paper we have introduced an FE model of curved membranes using higher dimensional shape functions that are
restricted to (the approximation of) the membrane surface. This allows for rapid insertion of arbitrarily shaped membranes into 
already existing 3D FE models, to be used for example for optimization purposes. We have shown numerically that the cut element approach to membranes
gives errors comparable to triangulated membranes, using the same degree of approximation, and we have proposed a stabilization method which
provides stability to the solution as well as giving the right conditioning of the discrete system, allowing for arbitrarily small cuts in the 3D mesh.

\section{Acknowledgements}
This research was supported in part by the Swedish Foundation for Strategic Research Grant No.\ AM13-0029, the Swedish Research Council Grants Nos.\ 2011-4992 and 2013-4708, 
and the Swedish Research Programme Essence

\newpage
\bibliographystyle{plain}
\bibliography{Membrane}
\end{document}